\documentclass[12pt,reqno]{amsart}

\usepackage{times}
\usepackage[T1]{fontenc}
\usepackage{mathrsfs}
\usepackage{latexsym}
\usepackage[dvips]{graphics}
\usepackage{epsfig}
\usepackage{amsmath,amsfonts,amsthm,amssymb,amscd}
\input amssym.def
\input amssym.tex
\usepackage{color}
\usepackage{verbatim}
\usepackage{hyperref}
\usepackage{pgf,tikz}
\usetikzlibrary{arrows}
\usepackage{subcaption}

\newcommand\be{\begin{equation}}
\newcommand\ee{\end{equation}}
\newcommand\bea{\begin{eqnarray}}
\newcommand\eea{\end{eqnarray}}
\newcommand\bi{\begin{itemize}}
\newcommand\ei{\end{itemize}}
\newcommand\ben{\begin{enumerate}}
\newcommand\een{\end{enumerate}}

\newtheorem{thm}{Theorem}[section]
\newtheorem{conj}[thm]{Conjecture}
\newtheorem{cor}[thm]{Corollary}
\newtheorem{lem}[thm]{Lemma}

\newtheorem{defi}[thm]{Definition}

\newtheorem{rek}[thm]{Remark}

\newcommand{\R}{\ensuremath{\mathbb{R}}}

\newcommand{\tri}{\bigtriangleup}
\newcommand{\nc}{\not \cong}
\newcommand{\ceiling}[1]{\left\lceil #1 \right\rceil}
\newcommand{\floor}[1]{\left\lfloor #1 \right\rfloor}

\numberwithin{equation}{section}

\begin{document}
\definecolor{qqqqff}{rgb}{0.,0.,1.}
\definecolor{ffqqqq}{rgb}{1.,0.,0.}

\title{Optimal point sets determining few distinct triangles}

\author{Alyssa Epstein}
\email{\textcolor{blue}{\href{mailto:ale2@williams.edu}{ale2@williams.edu}}} 
\address{Department of
  Mathematics and Statistics, Williams College, Williamstown, MA 01267}
  
\author{Adam Lott}
\email{\textcolor{blue}{\href{mailto:alott@u.rochester.edu} {alott@u.rochester.edu}}}
\address{Department of Mathematics, University of Rochester, Rochester, NY 14627}

\author{Steven J. Miller}
\email{\textcolor{blue}{\href{mailto:sjm1@williams.edu} {sjm1@williams.edu}, \href{mailto:Steven.Miller.MC.96@aya.yale.edu} {Steven.Miller.MC.96@aya.yale.edu}}}
\address{Department of
  Mathematics and Statistics, Williams College, Williamstown, MA 01267}
  
\author{Eyvindur A. Palsson}
\email{\textcolor{blue}{\href{mailto:eap2@williams.edu} {palsson@vt.edu}}}
\address{Department of Mathematics, Virginia Tech University, Blacksburg, VA 24061}

%\author{Steven J. Miller\thanks{E-mail: \texttt{sjmiller@math.brown.edu}}}



\date{\today}

\thanks{This work was supported by NSF Grants DMS1265673, DMS1561945, and DMS1347804, Williams College,  and the Clare Boothe Luce program. We also thank Paul Baird-Smith and Xiaoyu Xu for helpful conversations.}

\begin{abstract}

We generalize work of Erd\H{o}s and Fishburn to study the structure of finite point sets that determine few distinct triangles.  Specifically, we ask for a given $t$, what is the maximum number of points that can be placed in the plane to determine exactly $t$ distinct triangles?  Denoting this quantity by $F(t)$, we show that $F(1) = 4$, $F(2) = 5$, and 
%$F(t) < 48(t+1)$ for all $t$.  We 
we completely characterize the optimal configurations for $t = 1, 2$.  We also discuss the general structure of optimal configurations and conjecture that regular polygons are always optimal.  This differs from the structure of optimal configurations for distances, where it is conjectured that optimal configurations always exist in the triangular lattice.

\end{abstract}

\maketitle

\tableofcontents

\section{Introduction}

\begin{figure}[h]
\includegraphics[scale=0.8]{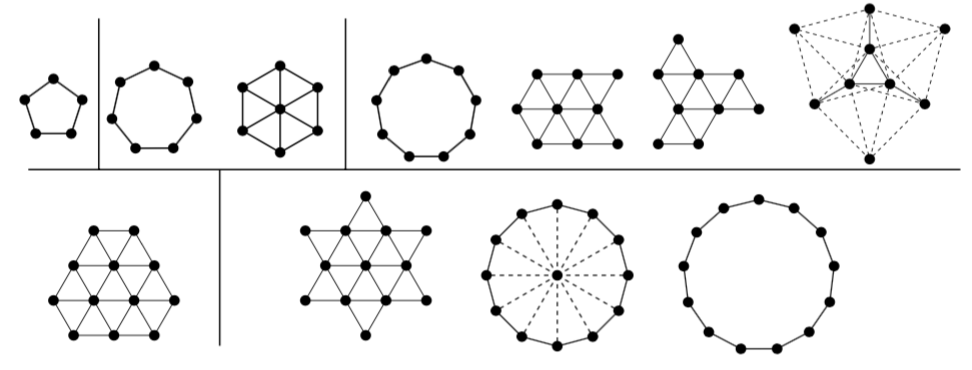}
\caption{Maximal configurations determining exactly $k$ distances, for $2 \leq k \leq 6$ \cite{BMP}.  For each $k>2$, there is an example from the triangular lattice; it is conjectured that this is always the case for $k$ large enough.}
\label{fig: distances}
\end{figure}

Finite point configurations are a central object of study in discrete geometry.  Perhaps the most well-known problem is the Erd\H{o}s distinct distances conjecture, which states that any set of $n$ points in the plane determines at least $\Omega(n/\sqrt{\log n})$ distinct distances between points.  This problem, first proposed by Erd\H{o}s in 1946 \cite{Erdos}, was essentially resolved by Guth and Katz who proved that $n$ points determined at least $\Omega(n/\log n)$ distinct distances \cite{GuthKatz}. Higher dimensional analogs still remain open. A closely related question is: given a fixed positive integer $k$, what is the maximum number of points that can be placed in the plane to determine exactly $k$ distances?  Furthermore, can the optimal configurations be completely characterized?  Erd\H{o}s and Fishburn \cite{ErdosFishburn} introduced this question in 1996 and characterized the optimal configurations for $1 \leq k \leq 4$.  Shinohara \cite{Sh} and Wei \cite{We} have characterized the optimal configurations for $k=5$ and $k=6$, respectively.  Erd\H{o}s also conjectured that an optimal configuration always exists in the triangular lattice given $k$ large enough (see Figure \ref{fig: distances}) and this conjecture remains open.

As a distance is just a pair of points, distances can be phrased as the set of 2-point configurations determined by a set.  Analogously, we can study the set of 3-point configurations (i.e., triangles) determined by a set.  %Much less is known about distinct triangles determined by a finite set.  
The analogue of the Erd\H{o}s distinct distance problem would ask for the minimum number of distinct triangles determined by n points in the plane.  It follows from Guth and Katz's result on the number of distinct distances that a set of $n$ points in the plane determines at least $\Omega(n^2)$ distinct triangles (see, for example, \cite{Ru}).  It is also known that this bound is best possible up to the implicit constant.  We study the following analogue of Erd\H{o}s and Fishburn's question: given a fixed $t$, what is the maximum number of points that can be placed in the plane to determine exactly $t$ distinct triangles?  Our main result is the following.

%In general, the current best known bound is that $n$ points determine at least $\Omega(n^{5/3})$ distinct triangles, a fact that follows from the best-known bound for the number of unit distances determined by a point set in the plane (see, for example, \cite{Szekely}).  Greenleaf and Iosevich \cite{GreenleafIosevich} show that this can be improved to $\Omega(n^{12/7 - \epsilon})$ if some restrictions are imposed on the set $P$.  

%In this paper, we adopt a slightly different definition of triangle than that of \cite{GreenleafIosevich}.  Loosely (this is made precise in the next section), we define a triangle as a triple of noncollinear points, and two triangles are considered distinct if and only if they are not congruent.  In \cite{GreenleafIosevich}, degenerate triangles (i.e., collinear triples of points) are allowed, but here they are not.

\begin{thm}
\label{thm: MainResult}
Let $F(t)$ denote the maximum number of points that can be placed in the plane to determine exactly $t$ distinct triangles. Then
\begin{enumerate}
\item \label{thmpart: 1Triangle}
$F(1) = 4$ and the only configuration that achieves this is a rectangle, and
\item \label{thmpart: 2Triangles}
$F(2) = 5$ and the only configurations that achieve this are a square with its center and a regular pentagon.
%%%%%%%%%%%%%%%%%%%%%%%%%%%%%%%%
\begin{comment}
\item \label{thmpart: UpperBound}
$F(t) < 48(t+1)$ for all $t$.
\end{comment}
%%%%%%%%%%%%%%%%%%%%%%%%%%%%%%%%%
\end{enumerate}
\end{thm}

We also make two conjectures: first, that $F(3) = 6$, with a regular hexagon being a representative optimal configuration, and second, that a regular polygon always minimizes the number of distinct triangles in an $n$-point set. If true, this second conjecture determines the true leading constant for Guth and Katz's asymptotic of at least $\Omega(n^2)$ distinct triangles for a set of $n$ points: $1/12$. 

We prove Theorem \ref{thm: MainResult} by classifying all potential arrangements of 4-point sets in the plane and sorting them by the minimum number of distinct triangles they create. To show part \ref{thmpart: 1Triangle}, we look at the 4-point sets that do not trivially determine more than one triangle. Through elementary geometry, we eliminate all non-trivial cases that have at least two distinct triangles except the rectangle. This immediately implies that $F(1)=4$, and the rectangle uniquely satisfies this equation. Proving part \ref{thmpart: 2Triangles}, we take the 4-point sets that determine fewer than three distinct triangles, and we examine all possible ways to add a fifth point to the set. After removing all cases where the fifth point causes at least three distinct triangles, the only remaining configurations are the square with a point at its center and the regular pentagon. Thus, $F(2)=5$.

\section{Conjectures}
In this section, we present some conjectures and investigate their consequences.

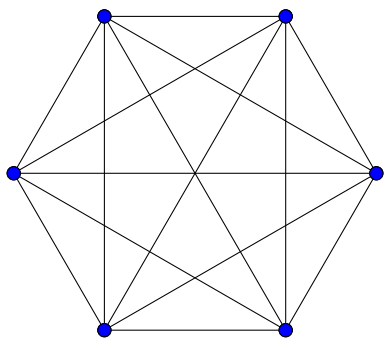
\begin{figure}[h]
\begin{tikzpicture}[line cap=round,line join=round,>=triangle 45,x=0.5cm,y=0.5cm]
\clip(-1.3437011700913217,-6.619671463395507) rectangle (10.240537924353422,3.547141267991481);
\draw (1.72,-5.34)-- (6.540663854698855,-5.34);
\draw (8.950995782048281,-1.1651826387253745)-- (6.540663854698855,-5.34);
\draw (8.950995782048281,-1.1651826387253745)-- (6.540663854698851,3.009634722549248);
\draw (1.72,3.0096347225492437)-- (6.540663854698851,3.009634722549248);
\draw (1.72,3.0096347225492437)-- (-0.6903319273494275,-1.1651826387253834);
\draw (1.72,-5.34)-- (-0.6903319273494275,-1.1651826387253834);
\draw (1.72,3.0096347225492437)-- (8.950995782048281,-1.1651826387253745);
\draw (1.72,3.0096347225492437)-- (6.540663854698855,-5.34);
\draw (1.72,3.0096347225492437)-- (1.72,-5.34);
\draw (-0.6903319273494275,-1.1651826387253834)-- (6.540663854698851,3.009634722549248);
\draw (-0.6903319273494275,-1.1651826387253834)-- (8.950995782048281,-1.1651826387253745);
\draw (-0.6903319273494275,-1.1651826387253834)-- (6.540663854698855,-5.34);
\draw (1.72,-5.34)-- (6.540663854698851,3.009634722549248);
\draw (1.72,-5.34)-- (8.950995782048281,-1.1651826387253745);
\draw (6.540663854698855,-5.34)-- (6.540663854698851,3.009634722549248);
\begin{scriptsize}
\draw [fill=qqqqff] (1.72,-5.34) circle (2.5pt);
\draw [fill=qqqqff] (6.540663854698855,-5.34) circle (2.5pt);
\draw [fill=qqqqff] (8.950995782048281,-1.1651826387253745) circle (2.5pt);
\draw [fill=qqqqff] (6.540663854698855,-5.34) circle (2.5pt);
\draw [fill=qqqqff] (8.950995782048281,-1.1651826387253745) circle (2.5pt);
\draw [fill=qqqqff] (6.540663854698851,3.009634722549248) circle (2.5pt);
\draw [fill=qqqqff] (1.72,3.0096347225492437) circle (2.5pt);
\draw [fill=qqqqff] (6.540663854698851,3.009634722549248) circle (2.5pt);
\draw [fill=qqqqff] (1.72,3.0096347225492437) circle (2.5pt);
\draw [fill=qqqqff] (-0.6903319273494275,-1.1651826387253834) circle (2.5pt);
\draw [fill=qqqqff] (1.72,-5.34) circle (2.5pt);
\draw [fill=qqqqff] (-0.6903319273494275,-1.1651826387253834) circle (2.5pt);
\end{scriptsize}
\end{tikzpicture}
\caption{A regular hexagon determines three distinct triangles.}
\label{fig: RegularHexagon}
\end{figure}

\begin{conj} \label{conj: 3TriangleSets}
Any set of seven points in the plane determines at least four distinct triangles; thus $F(3) = 6$.
\end{conj}
In Figure \ref{fig: RegularHexagon} we see that the vertices of a regular hexagon determine exactly three distinct triangles, so we know $F(3) \geq 6$. 

Another interesting question to ask concerns the general structure of the optimal configurations.  For example, are regular polygons always optimal?  What about regular polygons with their centers?  As we discussed in the introduction, Erd\H{o}s and Fishburn conjectured in \cite{ErdosFishburn} that optimal configurations for distinct distances always exist in the triangular lattice.  For triangles, we make an analogous but qualitatively different conjecture.

\begin{conj} \label{conj: RegularPolygons} \text{}
The regular $n$-gon minimizes (not necessarily uniquely) the number of distinct triangles determined by an $n$-point set.
\end{conj}

If true, Conjecture \ref{conj: RegularPolygons} establishes the following best-possible result on the number of distinct triangles which we prove in Section \ref{sec: TheoremProofDistinctTrianglesLowerBound}.

\begin{thm}
\label{thm: DistinctTrianglesLowerBound}
Unconditionally, the vertices of a regular $n$-gon determine $[n^2/12]$ distinct triangles, where $[y]$ denotes the nearest integer to $y$.  Assuming Conjecture \ref{conj: RegularPolygons}, this implies that $[n^2/12]$ is the minimum number of distinct triangles that can be determined by a set of $n$ points in the plane.
\end{thm}

\begin{rek}
It is known from the work of Guth and Katz that a set of $n$ points in the plane determines at least $\Omega(n^2)$ distinct triangles, and that this bound is best possible.  If true, Conjecture \ref{conj: RegularPolygons} establishes the true leading constant, namely $1/12$.
\end{rek}

\section{Definitions and setup}

We make precise the notion of distinct triangles.

\begin{defi}
Given a finite point set $P \subset \R^2$, we say two triples $(a,b,c), (a',b',c') \in P^3$ are equivalent if there is an isometry mapping one to the other, and we denote this as $(a,b,c) \sim (a',b',c')$.
\end{defi}

\begin{defi}
Given a finite point set $P \subset \R^2$, we denote by $P_{nc}^3$ the set of noncollinear triples $(a,b,c) \in P^3$.
\end{defi}

\begin{defi}
Given a finite point set $P \subset \R^2$, we define the set of distinct triangles determined by $P$ as
\begin{equation}
T(P) := P_{nc}^3 / \sim.
\end{equation}
\end{defi}

We prove %the first two parts of 
Theorem \ref{thm: MainResult} by enumerating cases and disposing of them one by one via elementary geometry.  
%We prove the last part by providing a self-contained lower bound for the minimum number of distinct triangles.  
We then conclude with a conjecture analogous to that of Erd\H{o}s concerning the structure of optimal configurations in general.

In the proof of Theorem \ref{thm: MainResult}, we also use the following lemma, which we prove in Section \ref{sec: LemmaProof}.

\begin{lem} \label{lem: FourPointConfigs}
For a set of four noncollinear points in the plane, exactly one of the following holds.
\begin{enumerate}
\item \label{part: NotConvex}
The four points are not in convex position.
\item 
The four points are in convex position.
\begin{enumerate}
\item \label{part: Collinear}
Three of the points are collinear.
\item \label{part: AllDistinct}
The determined quadrilateral has four distinct side lengths.
\item
The determined quadrilateral has exactly one pair of congruent sides.
\begin{enumerate}
\item \label{part: OnePairAdjacent}
The congruent sides are adjacent.
\item \label{part: OnePairOpposite}
The congruent sides are opposite.
\end{enumerate}
\item The determined quadrilateral has two distinct pairs of congruent sides.
\begin{enumerate}
\item \label{part: TwoPairAdjacent}
The congruent sides are adjacent to each other (a kite).
\item \label{part: TwoPairOpposite}
The congruent sides are opposite each other (a parallelogram).
\end{enumerate}
\item \label{part: ThreeSides}
Three sides are congruent and the fourth is distinct.
\item \label{part: FourSides}
All four sides are congruent (a rhombus).
\end{enumerate}
\end{enumerate}
Cases \ref{part: AllDistinct}, \ref{part: OnePairAdjacent}, \ref{part: OnePairOpposite}, and \ref{part: TwoPairAdjacent} determine at least three distinct triangles.  Cases \ref{part: NotConvex}, \ref{part: Collinear}, and \ref{part: ThreeSides} determine at least two distinct triangles.
\end{lem}

\section{Classifying optimal 1-triangle sets}

In this section, we prove part (\ref{thmpart: 1Triangle}) of Theorem \ref{thm: MainResult}.  We show that the only four-point configuration that determines exactly one triangle is a rectangle.  This proves that $F(1) = 4$ because there is no five-point configuration such that every four-point subconfiguration is a rectangle.

By Lemma \ref{lem: FourPointConfigs}, we only need to consider the cases \ref{part: TwoPairOpposite} and \ref{part: FourSides} because all of the other cases trivially lead to at least two triangles. We consider first the case \ref{part: TwoPairOpposite}, when there are two pairs of congruent sides opposite each other. 

\begin{proof}[Proof of case \ref{part: TwoPairOpposite}: two pairs of opposite congruent sides.]

\definecolor{qqwuqq}{rgb}{0.,0.39215686274509803,0.}
\definecolor{qqqqff}{rgb}{0.,0.,1.}

Since two pairs of opposite sides are congruent, the quadrilateral must be a parallelogram (Figure \ref{fig:parallelogram}).  We claim $\tri ABC$ and $\tri BCD$ are congruent if and only if $ABCD$ is a rectangle.  They share side $BC$ and $AB = CD$, so $\tri ABC \cong \tri BCD$ if and only if $BD = AC$, which happens if and only if $ABCD$ is a rectangle.
\end{proof}

\begin{comment}
When a parallelogram is not rectangular, it must form two distinct triangles: $\tri ABC$ and $\tri BCD$. The parallelogram has two defining angles, one obtuse, the other acute. In the image, these are $\angle ABC$ and $\angle BCD$, respectively. To prove $\tri ABC$ and $\tri BCD$ are distinct, we will show that either $\tri BCD$ is acute or right, or that the obtuse angle in it is smaller than the obtuse angle in $\tri ABC$. Without loss of generality, we focus on $\angle CBD$, though the argument can be repeated for $\angle BDC$, which we assume is acute for all of the following statements (we know that it cannot be obtuse or right at the same time as the other angle because they're part of the same triangle). If we let $\angle CBD$ be acute or right, we are done, because $\tri BCD$ has no obtuse angles. If $\angle CBD$ is obtuse, since it is contained within $\angle ABC$, $m\angle ABC > m\angle CBD$. (If $\angle ABC$ and $\angle CBD$ were the same, we would have a triangle, not a parallelogram). Because there is at most one obtuse angle in a triangle, the fact that these two obtuse angles are different in measure proves the existence of two distinct triangles. 

If, on the other hand, the parallelogram is a rectangle, it determines precisely one distinct triangle. 
\end{comment}

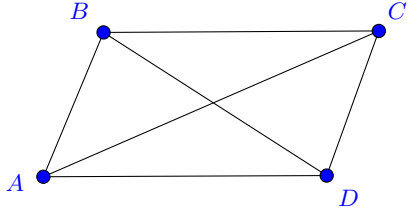
\begin{figure}[h]
\begin{tikzpicture}[line cap=round,line join=round,>=triangle 45,x=1.0cm,y=1.0cm]
\clip(-0.5,0) rectangle (5.489710896109946,3);
%\draw [shift={(1.,2.34)},color=qqwuqq,fill=qqwuqq,fill opacity=0.1] (0,0) -- (-112.61986494804043:0.30789487093842366) arc (-112.61986494804043:0.31308857502066406:0.30789487093842366) -- cycle;
%\draw [shift={(4.66,2.36)},color=qqwuqq,fill=qqwuqq,fill opacity=0.1] (0,0) -- (-179.68691142497934:0.30789487093842366) arc (-179.68691142497934:-109.81978157203878:0.30789487093842366) -- cycle;
\draw (0.2,0.42)-- (1.,2.34);
\draw (1.,2.34)-- (4.66,2.36);
\draw (4.66,2.36)-- (3.966437070439657,0.4356397288289089);
\draw (0.2,0.42)-- (3.966437070439657,0.4356397288289089);
\draw (1.,2.34)-- (3.966437070439657,0.4356397288289089);
\draw (0.2,0.42)-- (4.66,2.36);
\begin{scriptsize}
\draw [fill=qqqqff] (0.2,0.42) circle (2.5pt);
\draw[color=qqqqff] (-0.17555472915704978,0.3362312470823977) node {$A$};
\draw [fill=qqqqff] (1.,2.34) circle (2.5pt);
\draw[color=qqqqff] (0.6762877471059223,2.6249164543913435) node {$B$};
\draw [fill=qqqqff] (4.66,2.36) circle (2.5pt);
\draw[color=qqqqff] (4.904710641326941,2.6249164543913435) node {$C$};
\draw [fill=qqqqff] (3.966437070439657,0.4356397288289089) circle (2.5pt);
\draw[color=qqqqff] (4.258131412356251,0.1412311621547297) node {$D$};
\end{scriptsize}
\end{tikzpicture}
\caption{A quadrilateral with two pairs of opposite congruent sides.  If $ABCD$ is a rectangle, then it determines only one triangle, but if $ABCD$ is not a rectangle, then $\tri ABC$ and $\tri BCD$ are distinct. }
\label{fig:parallelogram}
\end{figure}

\begin{comment}
\begin{figure}[h]

\definecolor{qqqqff}{rgb}{0.,0.,1.}
\begin{tikzpicture}[line cap=round,line join=round,>=triangle 45,x=1.0cm,y=1.0cm]
\clip(-1.1483693634545042,-0.7546963409055187) rectangle (1.719238509603337,1.5499750872735154);
\draw (-0.6829958098022084,0.9234051214800301)-- (1.183786000291751,0.9440325447959854);
\draw (-0.6829958098022084,0.9234051214800301)-- (-0.6726820981442307,-0.025456351053914163);
\draw (-0.6726820981442307,-0.025456351053914163)-- (1.1940997119497287,-0.004828927737958852);
\draw (1.183786000291751,0.9440325447959854)-- (1.1940997119497287,-0.004828927737958852);
\begin{scriptsize}
\draw [fill=qqqqff] (-0.6829958098022084,0.9234051214800301) circle (2.5pt);
\draw[color=qqqqff] (-0.830815471471589,1.1169470527513587) node {$A$};
\draw [fill=qqqqff] (1.183786000291751,0.9440325447959854) circle (2.5pt);
\draw[color=qqqqff] (1.2573419394463692,1.150627010991971) node {$B$};
\draw [fill=qqqqff] (-0.6726820981442307,-0.025456351053914163) circle (2.5pt);
\draw[color=qqqqff] (-0.8452497392889942,-0.12439997954549062) node {$C$};
\draw [fill=qqqqff] (1.1940997119497287,-0.004828927737958852) circle (2.5pt);
\draw[color=qqqqff] (1.3054561655043866,-0.11477713433388712) node {$D$};
\end{scriptsize}
\end{tikzpicture}
\caption{Rectangle}
\label{fig:rectangle}
\end{figure}
\end{comment}

\begin{comment}
Since $AB$ = $CD$ and $AC$ = $BD$, and $\angle ABC$ = $\angle BDC$ = $\angle DCB$ = $\angle CAB$, by side-angle-side, $\tri ABC \cong \tri BDC \cong \tri DCB \cong \tri CAB$.
\end{comment}

\begin{proof}[Proof of case \ref{part: FourSides}: four congruent sides.]

Any quadrilateral with four sides congruent is a rhombus, and a rhombus is a parallelogram. So, by the argument in case \ref{part: TwoPairOpposite}, a rhombus determines two distinct triangles if and only if it is not a square.  Thus, we have shown that the only four-point configuration that determines one triangle is a rectangle.  This completes the proof of part (\ref{thmpart: 1Triangle}) of Theorem \ref{thm: MainResult}. 
\end{proof}

\section{Classifying optimal 2-triangle sets}

In this section, we prove part (\ref{thmpart: 2Triangles}) of Theorem \ref{thm: MainResult}.  As in the proof of part (\ref{thmpart: 1Triangle}), we show that the only possible configurations determining exactly two triangles are the square with its center and the regular pentagon.  We consider the possible four-point configurations enumerated in Lemma \ref{lem: FourPointConfigs}, and we show that the addition of a fifth point to any of them (unless it creates one of the two claimed configurations) necessarily determines a third triangle.  Moreover, adding a sixth point to either of the demonstrated optimal configurations also must determine a third triangle. By Lemma \ref{lem: FourPointConfigs}, the only cases we need to consider are \ref{part: NotConvex}, \ref{part: Collinear}, \ref{part: TwoPairOpposite}, \ref{part: ThreeSides}, and \ref{part: FourSides} because the other four point configurations already contain more than two distinct triangles.

\begin{figure}[h]
\begin{tikzpicture}[line cap=round,line join=round,>=triangle 45,x=0.6cm,y=0.6cm]
\clip(6.873764087152547,3.48760330578512) rectangle (16.34033057851246,14.8);
\draw (10.88732050807569,13.919934640057527)-- (8.22,9.34);
\draw (8.22,9.34)-- (13.52,9.32);
\draw (13.52,9.32)-- (10.88732050807569,13.919934640057527);
\draw (10.88732050807569,13.919934640057527)-- (10.875773502691898,10.859978213352507);
\draw (10.875773502691898,10.859978213352507)-- (8.22,9.34);
\draw (10.875773502691898,10.859978213352507)-- (13.52,9.32);
\draw [dash pattern=on 2pt off 2pt] (10.864226497308103,7.80002178664749)-- (8.22,9.34);
\draw [dash pattern=on 2pt off 2pt] (10.864226497308103,7.80002178664749)-- (13.52,9.32);
\draw [dash pattern=on 2pt off 2pt] (8.22,9.34)-- (10.852679491924315,4.740065359942471);
\draw [dash pattern=on 2pt off 2pt] (10.852679491924315,4.740065359942471)-- (13.52,9.32);
\begin{scriptsize}
\draw [fill=qqqqff] (8.22,9.34) circle (2.5pt);
\draw[color=qqqqff] (7.339579263711526,9.347858752817423) node {$A$};
\draw [fill=qqqqff] (13.52,9.32) circle (2.5pt);
\draw[color=qqqqff] (14,9.272727272727264) node {$B$};
\draw [fill=qqqqff] (10.88732050807569,13.919934640057527) circle (2.5pt);
\draw[color=qqqqff] (10.885785123966986,14.5) node {$C$};
\draw [fill=qqqqff] (10.875773502691898,10.859978213352507) circle (2.5pt);
\draw[color=qqqqff] (11.276468820435808,11.120961682945143) node {$D$};
\draw [fill=qqqqff] (10.864226497308103,7.80002178664749) circle (2.5pt);
\draw[color=qqqqff] (10.900811419985017,7.349361382419226) node {$E'$};
\draw [fill=qqqqff] (10.852679491924315,4.740065359942471) circle (2.5pt);
\draw[color=qqqqff] (10.84070623591289,4.17881292261457) node {$E$};
\end{scriptsize}
\end{tikzpicture}
\caption{Possibilities for adding a fifth point to a non-convex set.  } \label{fig: AddingNonConvex}
\end{figure}
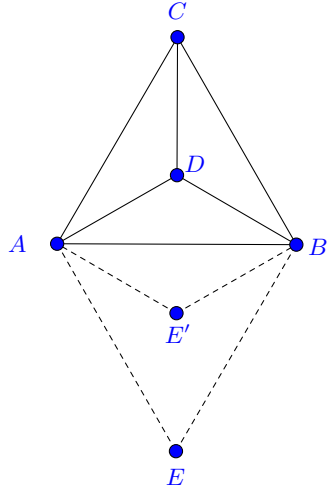

\begin{proof}[Proof of case \ref{part: NotConvex}: not in convex position.]

Using the notation of Figure \ref{fig: AddingNonConvex},  if $\tri ABC$ is not equilateral, or if $\tri ABC$ is equilateral but $D$ is not the center of $\tri ABC$, then there are already three distinct triangles, so no more work is needed.  

If $\tri ABC$ is equilateral and $D$ is its center, we show that the addition of a fifth point anywhere necessarily determines a new triangle.  When we add a fifth point $E$, it will necessarily determine a triangle with $AB$ (Figure \ref{fig: AddingNonConvex}).  If $\tri EAB$ is not congruent to $\tri ABC$ or $\tri ABD$, we're done, so assume it's congruent to one of those.  Either way, $\tri ECB$ will be distinct from the other two, so we have three distinct triangles, so this case is done. 
\end{proof}

\begin{figure}[h]
\begin{tikzpicture}[line cap=round,line join=round,>=triangle 45,x=0.7cm,y=0.7cm]
\clip(7.143567338586113,6.474598735414156) rectangle (14.589035886081072,11.7);
\draw (8.,9.)-- (13.822806055908455,9.);
\draw (8.,9.)-- (10.911403027954227,10.775356874530418);
\draw (10.911403027954227,10.775356874530418)-- (13.822806055908455,9.);
\draw (10.911403027954227,10.775356874530418)-- (10.911403027954226,9.);
\draw [dash pattern=on 2pt off 2pt] (8.,9.)-- (10.911403027954227,7.224643125469582);
\draw [dash pattern=on 2pt off 2pt] (10.911403027954227,7.224643125469582)-- (10.911403027954226,9.);
\draw [dash pattern=on 2pt off 2pt] (8.,9.)-- (8.,7.22464312546958);
\draw [dash pattern=on 2pt off 2pt] (8.,7.22464312546958)-- (10.911403027954226,9.);
\begin{scriptsize}
\draw [fill=qqqqff] (8.,9.) circle (2.5pt);
\draw[color=qqqqff] (7.624151531921529,9.078360857514987) node {$A$};
\draw [fill=qqqqff] (13.822806055908455,9.) circle (2.5pt);
\draw[color=qqqqff] (14.173007778417578,9.049669263883022) node {$B$};
\draw [fill=qqqqff] (10.911403027954226,9.) circle (2.5pt);
\draw[color=qqqqff] (11.3,9.3) node {$C$};
\draw [fill=qqqqff] (10.911403027954227,10.775356874530418) circle (2.5pt);
\draw[color=qqqqff] (10.887820307557568,11.180020091056429) node {$D$};
\draw [fill=qqqqff] (10.911403027954227,7.224643125469582) circle (2.5pt);
\draw[color=qqqqff] (10.923684799597524,6.8) node {$E$};
\draw [fill=qqqqff] (8.,7.22464312546958) circle (2.5pt);
\draw[color=qqqqff] (7.5,6.9) node {$E'$};
\draw[fill=qqqqff] (8,10.8) circle (2.5pt);
\draw[color=qqqqff] (8,11.2) node {$E''$};
\draw [dash pattern=on 2pt off 2pt] (8.,10.8)-- (10.911403027954226,9.);
\draw [dash pattern=on 2pt off 2pt] (8.,10.8)-- (8,9.);
\end{scriptsize}
\end{tikzpicture}
\caption{Addition of a fifth point when three points are collinear.  If $DC \neq AC$, then any choice of $E$ forces a third triangle.  If, on the other hand, $DC = AC$, then choosing $E$ creates a square with its center but $E'$ and $E''$ still generate a third triangles.} \label{fig: AddingCollinear}
\end{figure}
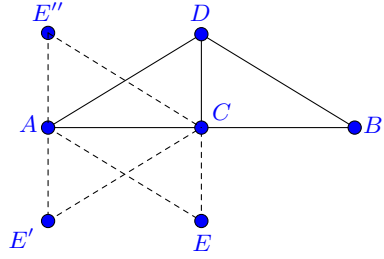

\begin{proof}[Proof of case \ref{part: Collinear}: three collinear points.]

With the notation of Figure \ref{fig: AddingCollinear}, if $D$ does not lie on the perpendicular bisector of $AB$, then $\tri ACD$, $\tri BCD$, and $\tri ABD$ are all distinct, so no more work is needed.  Also note that if a fifth point $E$ is added to the interior of $\tri  ABD$, it creates a non-convex four-point subconfiguration, so the previous case applies to show that there are at least 3 distinct triangles.  Thus we assume the fifth point $E$ is added  outside $\tri ABD$.

If $D$ lies on the perpendicular bisector of $AB$ but $DC \neq AB$, the addition of a fifth point $E$ will create a triangle with $AC$.  Triangle $\tri EAC$ can't be congruent to $\tri ABD$ because $AC$ is shorter than any side of $\tri ABD$, so to avoid a third triangle we must have $\tri EAC \cong \tri ACD$.  There are three choices for $E$ that satisfy this (Figure \ref{fig: AddingCollinear}), but either way, $\tri EAC$, $\tri EAB$, and $\tri EDB$ are all distinct.

If $D$ lies on the perpendicular bisector of $AB$ and $DC = AB$, then the same argument from above still applies; however, in this case, choosing $E$ to form the square $ADBE$ leaves us with only two triangles, but the other two choices for $E$ give us three (see Figure \ref{fig: AddingCollinear}), so this case is done.
\end{proof}

\begin{proof}[Proof of case \ref{part: TwoPairOpposite}: two pairs of opposite congruent sides.]

This case has two subcases.

\emph{Subcase A: non-rectangle}: Using the notation of Figure \ref{fig: AddingTwoPairOpposite}, if we add a fifth point $E$ on line $AB$, then we have five points with three collinear, so we have 3 distinct triangles by case \ref{fig: AddingCollinear}.  So assume $E$ does not lie on line $AB$.  Then $\tri EAB$ will be created.  If $\tri EAB$ is distinct from both $\tri ABC$ and $\tri ABD$, then we also have three distinct triangles, so assume otherwise.  The only ways this can happen are enumerated in Figure \ref{fig: AddingTwoPairOpposite}.  In Figure \ref{fig: AddingOption1TwoPairOpposite}, point $E$ creates three collinear points ($EAD$), point $E'$ creates a non-convex subconfiguration ($ACBE'$), and point $E''$ creates three collinear points ($CDE''$).  Thus in any case there will be three distinct triangles.  In Figure \ref{fig: AddingOption2TwoPairOpposite}, point $E'$ creates three collinear points ($CBE'$) and point $E''$ also creates three collinear points ($DE''C$).  Point $E$ creates a kite $ADBE$ if $AD \neq DB$, and if $AD = DB$, then $CBE$ must be collinear, so in this case also, we have three distinct triangles no matter what.  

\emph{Subcase B: non-square rectangle}: If the fifth point is added inside the rectangle, then we get either a non-convex configuration or a configuration with three collinear points (Figure \ref{fig: AddingRectangleInside}).  So assume that the fifth point is added outside the rectangle.  Using the notation of Figure \ref{fig: AddingRectangleOutside}, to add a fifth point $E$ without creating three distinct triangles there are three potential possibilities.
\begin{enumerate}
\item $\tri EAB \cong \tri ABC$.  In this case, we get three collinear points, so we have three triangles.
\item $\tri E'AD \cong \tri ABC$.  Here, $DCE$ are collinear, so we have three triangles.
\item $\tri E''DC \cong \tri E''CB \nc \tri ABC$.  In this case, $E''DAB$ will form a kite, so we have three triangles.
\end{enumerate}

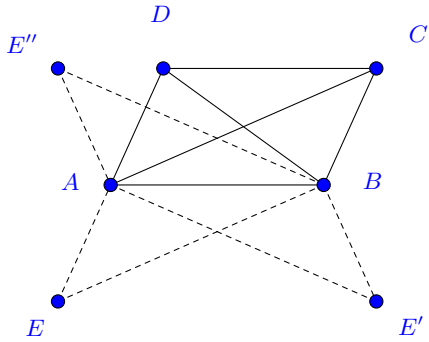
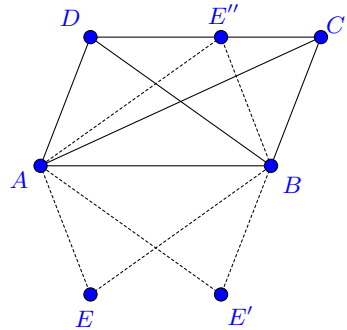
\begin{figure}[h]
\begin{subfigure}[h]{.4\textwidth}
\begin{tikzpicture}[line cap=round,line join=round,>=triangle 45,x=0.5cm,y=0.5cm]
\clip(3.9703531179564315,5.160916604057115) rectangle (16.862915101427504,16.280375657400466);
\draw (7.12,10.44)-- (12.784274004671742,10.44);
\draw (8.52,13.54)-- (14.184274004671742,13.54);
\draw (12.784274004671742,10.44)-- (14.184274004671742,13.54);
\draw (8.52,13.54)-- (7.12,10.44);
\draw (7.12,10.44)-- (14.184274004671742,13.54);
\draw (8.52,13.54)-- (12.784274004671742,10.44);
\draw [dash pattern=on 2pt off 2pt] (5.72,7.34)-- (7.12,10.44);
\draw [dash pattern=on 2pt off 2pt] (12.784274004671742,10.44)-- (14.184274004671742,7.34);
\draw [dash pattern=on 2pt off 2pt] (5.72,7.34)-- (12.784274004671742,10.44);
\draw [dash pattern=on 2pt off 2pt] (14.184274004671742,7.34)-- (7.12,10.44);
\draw [dash pattern=on 2pt off 2pt] (5.72,13.54)-- (7.12,10.44);
\draw [dash pattern=on 2pt off 2pt] (5.72,13.54)-- (12.784274004671742,10.44);
\begin{scriptsize}
\draw [fill=qqqqff] (7.12,10.44) circle (2.5pt);
\draw[color=qqqqff] (6.043981968444785,10.510277986476348) node {$A$};
\draw [fill=qqqqff] (12.784274004671742,10.44) circle (2.5pt);
\draw[color=qqqqff] (14.083050338091667,10.540330578512412) node {$B$};
\draw [fill=qqqqff] (8.52,13.54) circle (2.5pt);
\draw[color=qqqqff] (8.448189331329834,15.003140495867784) node {$D$};
\draw [fill=qqqqff] (14.184274004671742,13.54) circle (2.5pt);
\draw[color=qqqqff] (15.300180315552222,14.462193839218648) node {$C$};
\draw [fill=qqqqff] (5.72,7.34) circle (2.5pt);
\draw[color=qqqqff] (5.112351615326829,6.633493613824206) node {$E$};
\draw [fill=qqqqff] (14.184274004671742,7.34) circle (2.5pt);
\draw[color=qqqqff] (15.119864763335844,6.693598797896333) node {$E'$};
\draw [fill=qqqqff] (5.72,13.54) circle (2.5pt);
\draw[color=qqqqff] (4.766746806912104,14.19172051089408) node {$E''$};
\end{scriptsize}
\end{tikzpicture}
\caption{Possibilities for $E$ so that $\tri EAB \cong \tri ABC$.  Any one of these choices creates a 4-point subconfiguration determining at least 3 distinct triangles.}
\label{fig: AddingOption1TwoPairOpposite}
\end{subfigure}
\hspace{25mm}
\begin{subfigure}[h]{.4\textwidth}
\begin{tikzpicture}[line cap=round,line join=round,>=triangle 45,x=0.4cm,y=0.4cm]
\clip(5.124256539853833,2.150017075336393) rectangle (17.158953623386367,13.651963663684178);
\draw (6.24,8.06)-- (13.90065271370528,8.06);
\draw (7.9,12.34)-- (15.56065271370528,12.34);
\draw (6.24,8.06)-- (7.9,12.34);
\draw (13.90065271370528,8.06)-- (15.56065271370528,12.34);
\draw (6.24,8.06)-- (15.56065271370528,12.34);
\draw (7.9,12.34)-- (13.90065271370528,8.06);
\draw [dash pattern=on 1pt off 1pt] (6.24,8.06)-- (7.9,3.78);
\draw [dash pattern=on 1pt off 1pt] (7.9,3.78)-- (13.90065271370528,8.06);
\draw [dash pattern=on 1pt off 1pt] (13.90065271370528,8.06)-- (12.24065271370528,3.78);
\draw [dash pattern=on 1pt off 1pt] (12.24065271370528,3.78)-- (6.24,8.06);
\draw [dash pattern=on 1pt off 1pt] (13.90065271370528,8.06)-- (12.24065271370528,12.34);
\draw [dash pattern=on 1pt off 1pt] (12.24065271370528,12.34)-- (6.24,8.06);
\begin{scriptsize}
\draw [fill=qqqqff] (6.24,8.06) circle (2.5pt);
\draw[color=qqqqff] (5.520404343965573,7.586804180042355) node {$A$};
\draw [fill=qqqqff] (13.90065271370528,8.06) circle (2.5pt);
\draw[color=qqqqff] (14.604483300321004,7.409220681647437) node {$B$};
\draw [fill=qqqqff] (7.9,12.34) circle (2.5pt);
\draw[color=qqqqff] (7.227937982378248,12.968950208319107) node {$D$};
\draw [fill=qqqqff] (15.56065271370528,12.34) circle (2.5pt);
\draw[color=qqqqff] (16.05247182569495,12.764046171709587) node {$C$};
\draw [fill=qqqqff] (7.9,3.78) circle (2.5pt);
\draw[color=qqqqff] (7.706047401133797,2.9696332217744787) node {$E$};
\draw [fill=qqqqff] (12.24065271370528,3.78) circle (2.5pt);
\draw[color=qqqqff] (12.896949661908328,3.0652551055255888) node {$E'$};
\draw [fill=qqqqff] (12.24065271370528,12.34) circle (2.5pt);
\draw[color=qqqqff] (12.336878628508972,13.105552899392123) node {$E''$};
\end{scriptsize}
\end{tikzpicture}
\caption{Possibilities for $E$ so that $\tri EAB \cong \tri ABD$.  Here also, any choice creates a bad 4-point subconfiguration.}
\label{fig: AddingOption2TwoPairOpposite}
\end{subfigure}
\caption{Possible additions of a fifth point when two pairs of opposite sides are congruent.} \label{fig: AddingTwoPairOpposite}
\end{figure}

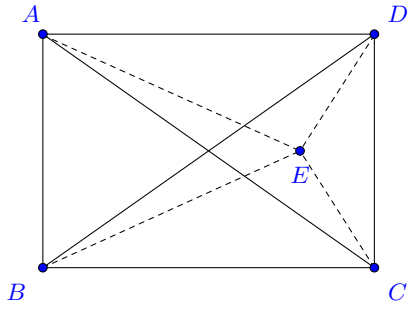
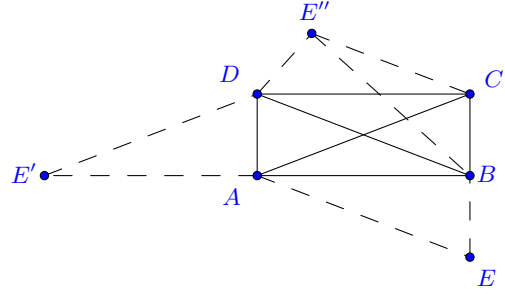
\begin{figure}[h]
%%%%%%%%%%%%%%%%%%%%%%%%%%
\begin{subfigure}{0.4\textwidth}
\begin{tikzpicture}[line cap=round,line join=round,>=triangle 45,x=1.0cm,y=1.0cm,scale=0.65]
\clip(-3,-1.8) rectangle (7.3,5.3);
%\fill[color=zzttqq,fill=zzttqq,fill opacity=0.1] (-1.26,4.06) -- (-1.24,-0.7) -- (3.52,-0.68) -- (3.5,4.08) -- cycle;
\draw (-1.26,4.06)-- (-1.26,-0.7);
\draw (-1.26,-0.7)-- (5.52,-0.7);
\draw (5.52,-0.7)-- (5.52,4.08);
\draw (5.52,4.08)-- (-1.26,4.08);
\draw (-1.26,4.08)-- (5.52,-0.7);
\draw (-1.26,-0.7)-- (5.52,4.08);
\draw [dash pattern=on 2pt off 2pt] (-1.26,4.06)-- (4,1.69);
\draw [dash pattern=on 2pt off 2pt] (-1.26,-0.7)-- (4,1.69);
\draw [dash pattern=on 2pt off 2pt] (5.52,4.06)-- (4,1.69);
\draw [dash pattern=on 2pt off 2pt] (5.52,-0.7)-- (4,1.69);
\begin{scriptsize}
\draw [fill=qqqqff] (-1.26,4.08) circle (2.5pt);
\draw[color=qqqqff] (-1.5,4.5) node {$A$};
\draw [fill=qqqqff] (-1.26,-0.7) circle (2.5pt);
\draw[color=qqqqff] (-1.8,-1.2) node {$B$};
\draw [fill=qqqqff] (5.52,-0.7) circle (2.5pt);
\draw [color=qqqqff] (6,-1.2) node {$C$};
\draw [fill=qqqqff] (5.52,4.08) circle (2.5pt);
\draw [color=qqqqff] (6,4.5) node {$D$};
\draw [fill=qqqqff] (4,1.69) circle (2.5pt);
\draw [color=qqqqff] (4,1.2) node {$E$};
\end{scriptsize}
\end{tikzpicture}
\caption{Any way to place a fifth point inside a rectangle results in at least 3 distinct triangles.}
\label{fig: AddingRectangleInside}
\end{subfigure}
%%%%%%%%%%%%%%%%%%%%%%%%%%%%
\hspace{15mm}
%%%%%%%%%%%%%%%%%%%%%%%%%%%%%%%%%%
\begin{subfigure}{0.4\textwidth}
\begin{tikzpicture}[line cap=round,line join=round,>=triangle 45,x=1.0cm,y=1.0cm,scale=0.35]
\clip(-10.981487981900603,-11.528183420236466) rectangle (9.80322011001988,0.5725027976624583);
\draw (-0.96,-6.3)-- (7.120891039978207,-6.3);
\draw (-0.96,-3.2)-- (-0.96,-6.3);
\draw (-0.96,-3.2)-- (7.120891039978207,-3.2);
\draw (7.120891039978207,-3.2)-- (7.120891039978207,-6.3);
\draw (-0.96,-3.2)-- (7.120891039978207,-6.3);
\draw (-0.96,-6.3)-- (7.120891039978207,-3.2);
\draw [dash pattern=on 7pt off 7pt] (7.120891039978207,-9.4)-- (7.120891039978207,-6.3);
\draw [dash pattern=on 7pt off 7pt] (-0.96,-6.3)-- (7.120891039978207,-9.4);
\draw [dash pattern=on 7pt off 7pt] (-0.96,-6.3)-- (-9.040891039978208,-6.3);
\draw [dash pattern=on 7pt off 7pt] (-9.040891039978208,-6.3)-- (-0.96,-3.2);
\draw [dash pattern=on 7pt off 7pt] (-0.96,-3.2)-- (1.113328889671195,-0.8953726298477651);
\draw [dash pattern=on 7pt off 7pt] (1.113328889671195,-0.8953726298477651)-- (7.120891039978207,-6.3);
\draw [dash pattern=on 7pt off 7pt] (1.113328889671195,-0.8953726298477651)-- (7.120891039978207,-3.2);
\begin{scriptsize}
\draw [fill=qqqqff] (-0.96,-6.3) circle (4.5pt);
\draw[color=qqqqff] (-1.927327470625652,-7.114991976061564) node {$A$};
\draw [fill=qqqqff] (7.120891039978207,-6.3) circle (4.5pt);
\draw[color=qqqqff] (7.753221503693477,-6.232353687226583) node {$B$};
\draw [fill=qqqqff] (-0.96,-3.2) circle (4.5pt);
\draw[color=qqqqff] (-2.0127440792225855,-2.417078503230217) node {$D$};
\draw [fill=qqqqff] (7.120891039978207,-3.2) circle (4.5pt);
\draw[color=qqqqff] (8.037943532349923,-2.616383923289729) node {$C$};
\draw [fill=qqqqff] (7.120891039978207,-9.4) circle (4.5pt);
\draw[color=qqqqff] (7.753221503693477,-10.189989885551173) node {$E$};
\draw [fill=qqqqff] (-9.040891039978208,-6.3) circle (4.5pt);
\draw[color=qqqqff] (-9.81412766440918,-6.203881484360939) node {$E'$};
\draw [fill=qqqqff] (1.113328889671195,-0.8953726298477651) circle (4.5pt);
\draw[color=qqqqff] (1.2330870474608873,-0.025413462516076818) node {$E''$};
\end{scriptsize}
\end{tikzpicture}
\caption{Any way to place a fifth point outside a rectangle also results in at least 3 distinct triangles.}
\label{fig: AddingRectangleOutside}
\end{subfigure}
%%%%%%%%%%%%%%%%%%
\caption{Any way to add a fifth point to a rectangle results in at least 3 distinct triangles.}
\label{fig: AddingRectangle}
\end{figure}
So we see both subcases yield at least three triangles, so the proof of case \ref{part: TwoPairOpposite} is complete.
\end{proof}

\begin{proof}[Proof of case \ref{part: ThreeSides}: three congruent sides.]
Using the notation of Figure \ref{fig: AddingThreeSidesCongruent}, if the quadrilateral $ABCD$ is not a trapezoid, then in particular $AC \neq BD$.  Then we claim $\tri ABD$, $\tri BDC$, and $\tri ABC$ are all distinct.  Triangle $\tri ABC \nc \tri ABD$ because $AC \neq BD$.  If $\tri ABC \cong \tri BDC$, then $AB = BD$ and $CD = AC$, but this is impossible because then there would be two isoceles triangles based on $AD$.

So we can assume $ABCD$ is a trapezoid.  When we add a fifth point $E$, $\tri EAD$ is created (Figure \ref{fig: AddingThreeSidesCongruent}).  As in case \ref{part: TwoPairOpposite}, we must have $\tri EAD \cong \tri ABD$ or $\tri EAD \cong \tri ACD$.  Suppose $\tri EAD \cong \tri ABD$ (Figure \ref{fig: Option1AddingThreeSidesCongruent}).  In the figure, point $E$ creates a non-convex configuration $EABD$ and point $E'$ creates three collinear points $E'DC$.  For point $E''$, if $E''C$ is a new distance then we obviously have a new triangle.  If $E''C = DC$, then $E''DC$ is a new triangle.  If $E''C = AC$, then $E''DAC$ is a kite, so we have three triangles.  If $E''C = BC$, then $ABCE''D$ is a regular pentagon, and this is one of our claimed optimal configurations.

Now suppose that $\tri EAD \cong \tri ACD$ (Figure \ref{fig: Option2AddingThreeSidesCongruent}).  Point $E$ in the figure makes $EACD$ either a kite, a non-convex congfiguration, or a configuration with three collinear points, depending on the length of $DC$.  In any case, we have at least three triangles.  Point $E'$ makes three collinear points $E'AB$.  For point $E''$, if $E''C$ is a new distance, we have a new triangle.  If $E''C = AD$, then $ADE''C$ is a non-rhombus parallelogram, so we have three triangles.  If $E''C = AC$, then $DE''C$ is a new triangle.  Finally, if $E''C = DC$, then $DE''C$ is also a new triangle.  This shows that the only way to add a fifth point to a trapezoid configuration without generating a third triangle is to create a regular pentagon, which concludes the proof of case \ref{part: ThreeSides}.  
\end{proof}

\begin{figure}[h]
\begin{subfigure}[h]{.4\textwidth}
\begin{tikzpicture}[line cap=round,line join=round,>=triangle 45,x=0.25cm,y=0.25cm,scale=0.8]
\clip(-6.3635864496637415,0.00775494018701372) rectangle (25.835836104285338,25.619404194998314);
\draw (9.02,8.72)-- (16.903653975156438,8.72);
\draw (9.02,8.72)-- (5.078173012421782,15.54744461713165);
\draw (20.845480962734655,15.54744461713165)-- (16.903653975156438,8.72);
\draw (5.078173012421782,15.54744461713165)-- (20.845480962734655,15.54744461713165);
\draw (5.078173012421782,15.54744461713165)-- (16.903653975156438,8.72);
\draw (9.02,8.72)-- (20.845480962734655,15.54744461713165);
\draw [dash pattern=on 2pt off 2pt] (9.02,8.72)-- (5.078173012421775,1.892555382868344);
\draw [dash pattern=on 2pt off 2pt] (5.078173012421775,1.892555382868344)-- (5.078173012421782,15.54744461713165);
\draw [dash pattern=on 2pt off 2pt] (5.078173012421772,15.547444617131646)-- (-2.805480962734667,15.547444617131651);
\draw [dash pattern=on 2pt off 2pt] (-2.805480962734667,15.547444617131651)-- (9.02,8.72);
\draw [dash pattern=on 2pt off 2pt] (5.078173012421777,15.547444617131648)-- (9.02,22.3748892342633);
\draw [dash pattern=on 2pt off 2pt] (9.02,22.3748892342633)-- (9.02,8.72);
\begin{scriptsize}
\draw [fill=qqqqff] (9.02,8.72) circle (2.5pt);
\draw[color=qqqqff] (10.173068974702588,6.797194768890272) node {$A$};
\draw [fill=qqqqff] (16.903653975156438,8.72) circle (2.5pt);
\draw[color=qqqqff] (18.34056342210303,6.898028033672993) node {$B$};
\draw [fill=qqqqff] (5.078173012421782,15.54744461713165) circle (2.5pt);
\draw[color=qqqqff] (3.2155737046948043,17.855242806728747) node {$D$};
\draw [fill=qqqqff] (20.845480962734655,15.54744461713165) circle (2.5pt);
\draw[color=qqqqff] (23.3486155729782,15.468855540204334) node {$C$};
\draw [fill=qqqqff] (5.078173012421775,1.892555382868344) circle (2.5pt);
\draw[color=qqqqff] (2.9130739103466397,1.284976294101488) node {$E$};
\draw [fill=qqqqff] (-2.805480962734667,15.547444617131651) circle (2.5pt);
\draw[color=qqqqff] (-5.086365095749268,15.569688804987056) node {$E'$};
\draw [fill=qqqqff] (9.02,22.3748892342633) circle (2.5pt);
\draw[color=qqqqff] (9.131125238614466,24.10690522325749) node {$E''$};
\end{scriptsize}
\end{tikzpicture}
\caption{Options for adding a fifth point $E$ so that $\tri EAD \cong \tri ABD$.  Adding $E$ or $E'$ will create a third triangle, and adding $E''$ will create a third triangle if and only if $DC \neq AC$.  If $DC = AC$, $E''$ is the fifth vertex of a regular pentagon.}
\label{fig: Option1AddingThreeSidesCongruent}
\end{subfigure}
\hspace{15mm}
\begin{subfigure}{.4\textwidth}
\begin{tikzpicture}[line cap=round,line join=round,>=triangle 45,x=0.25cm,y=0.25cm,scale=0.8]
\clip(-17.03681585299279,-6.4928188752347085) rectangle (18.497715104070608,17.451292278808804);
\draw (2.82,1.78)-- (10.703653975156438,1.78);
\draw (2.82,1.78)-- (-1.1218269875782174,8.60744461713165);
\draw (14.645480962734656,8.60744461713165)-- (10.703653975156438,1.78);
\draw (-1.1218269875782174,8.60744461713165)-- (14.645480962734656,8.60744461713165);
\draw (-1.1218269875782174,8.60744461713165)-- (10.703653975156438,1.78);
\draw (2.82,1.78)-- (14.645480962734656,8.60744461713165);
\draw [dash pattern=on 2pt off 2pt] (-1.1218269875782187,8.60744461713165)-- (-9.00548096273466,-5.047444617131643);
\draw [dash pattern=on 2pt off 2pt] (-9.00548096273466,-5.047444617131643)-- (2.82,1.78);
\draw [dash pattern=on 2pt off 2pt] (-12.947307950312874,1.78)-- (-1.1218269875782187,8.607444617131648);
\draw [dash pattern=on 2pt off 2pt] (-12.947307950312874,1.78)-- (2.82,1.78);
\draw [dash pattern=on 2pt off 2pt] (10.703653975156447,15.434889234263292)-- (-1.1218269875782152,8.607444617131652);
\draw [dash pattern=on 2pt off 2pt] (2.82,1.78)-- (10.703653975156447,15.434889234263292);
\begin{scriptsize}
\draw [fill=qqqqff] (2.82,1.78) circle (2.5pt);
\draw[color=qqqqff] (2.751705761919205,0.2918171070208349) node {$A$};
\draw [fill=qqqqff] (10.703653975156438,1.78) circle (2.5pt);
\draw[color=qqqqff] (11.910964337964185,0.2635477904281035) node {$B$};
\draw [fill=qqqqff] (-1.1218269875782174,8.60744461713165) circle (2.5pt);
\draw[color=qqqqff] (-2.676003023885228,10.553579030182341) node {$D$};
\draw [fill=qqqqff] (14.645480962734656,8.60744461713165) circle (2.5pt);
\draw[color=qqqqff] (16.74501747532126,8.54645755209841) node {$C$};
\draw [fill=qqqqff] (-9.00548096273466,-5.047444617131643) circle (2.5pt);
\draw[color=qqqqff] (-10.591411669850027,-5.107622362190866) node {$E$};
\draw [fill=qqqqff] (-12.947307950312874,1.78) circle (2.5pt);
\draw[color=qqqqff] (-14.351230776683305,1.648744303471943) node {$E'$};
\draw [fill=qqqqff] (10.703653975156447,15.434889234263292) circle (2.5pt);
\draw[color=qqqqff] (11.826156388185991,15.981287815986773) node {$E''$};
\end{scriptsize}
\end{tikzpicture}
\caption{Options for adding a fifth point $E$ so that $\tri EAD \cong \tri ACD$.  $E$ and $E'$ both generate a third triangle, and $E''$ generates a third triangle if $D \neq AC$.  If $DC = AC$, then $E''$ and $C$ are the same point.}
\label{fig: Option2AddingThreeSidesCongruent}
\end{subfigure}
\caption{Possible additions of a fifth point when three sides are congruent.}
\label{fig: AddingThreeSidesCongruent}
\end{figure}

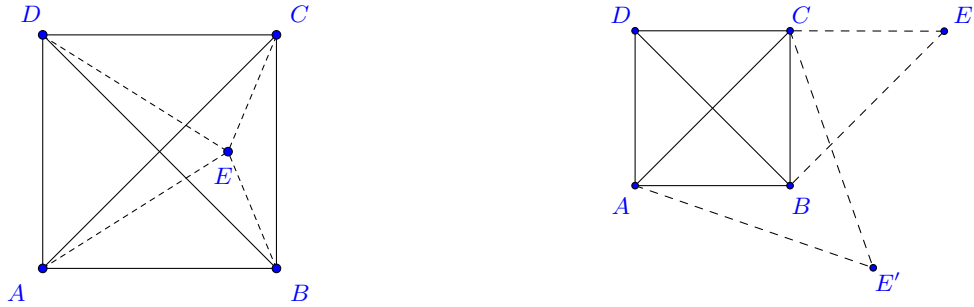
\begin{figure}[h]
%%%%%%%%%%%%%%%%%%%%%%%%%%%%%%%%%%%%%%%%%%%%%%%%%%%%%%%%%%%
\begin{subfigure}{0.4\textwidth}
\vspace{10mm}
\begin{tikzpicture}[line cap=round,line join=round,>=triangle 45,x=1.0cm,y=1.0cm,scale=0.65]
\clip(-3.3,-1.8) rectangle (7.3,4.8);
%\fill[color=zzttqq,fill=zzttqq,fill opacity=0.1] (-1.26,4.06) -- (-1.24,-0.7) -- (3.52,-0.68) -- (3.5,4.08) -- cycle;
\draw (-1.26,4.06)-- (-1.26,-0.7);
\draw (-1.26,-0.7)-- (3.52,-0.7);
\draw (3.52,-0.7)-- (3.52,4.08);
\draw (3.52,4.08)-- (-1.26,4.08);
\draw (-1.26,4.08)-- (3.52,-0.7);
\draw (-1.26,-0.7)-- (3.52,4.08);
\draw [dash pattern=on 2pt off 2pt] (-1.26,4.06)-- (2.53,1.69);
\draw [dash pattern=on 2pt off 2pt] (-1.26,-0.7)-- (2.53,1.69);
\draw [dash pattern=on 2pt off 2pt] (3.52,4.06)-- (2.53,1.69);
\draw [dash pattern=on 2pt off 2pt] (3.52,-0.7)-- (2.53,1.69);
\begin{scriptsize}
\draw [fill=qqqqff] (-1.26,4.08) circle (2.5pt);
\draw[color=qqqqff] (-1.5,4.5) node {$D$};
\draw [fill=qqqqff] (-1.26,-0.7) circle (2.5pt);
\draw[color=qqqqff] (-1.8,-1.2) node {$A$};
\draw [fill=qqqqff] (3.52,-0.7) circle (2.5pt);
\draw [color=qqqqff] (4,-1.2) node {$B$};
\draw [fill=qqqqff] (3.52,4.08) circle (2.5pt);
\draw [color=qqqqff] (4,4.5) node {$C$};
\draw [fill=qqqqff] (2.53,1.69) circle (2.5pt);
\draw [color=qqqqff] (2.43,1.2) node {$E$};
\end{scriptsize}
\end{tikzpicture}
\caption{Addition of a fifth point inside the square but not at the center.  $ABCE$ is a non-convex configuration, so we get three distinct triangles.}
\label{fig: AddingSquareInside}
\end{subfigure}
%%%%%%%%%%%%%%%%%%%%%%%%%%%%%%%%%%%%%%%%%%%%%%%%
\hspace{20mm}
%%%%%%%%%%%%%%%%%%%%%%%%%%%%%%%%%%%%%%%%%%%%%%
\begin{subfigure}{0.4\textwidth}
\definecolor{uuuuuu}{rgb}{0.26666666666666666,0.26666666666666666,0.26666666666666666}
\definecolor{zzttqq}{rgb}{0.6,0.2,0.}
\definecolor{qqqqff}{rgb}{0.,0.,1.}
\begin{tikzpicture}[line cap=round,line join=round,>=triangle 45,x=1.0cm,y=1.0cm,scale=0.5]
\clip(-6.8824065040650355,-3) rectangle (6.159173054587681,5.5);
%\fill[color=zzttqq,fill=zzttqq,fill opacity=0.1] (-4.623349593495929,0.0494680603948896) -- (-0.5033495934959349,0.0494680603948896) -- (-0.5033495934959347,4.169468060394882) -- (-4.623349593495928,4.169468060394884) -- cycle;
\draw (-4.623349593495929,0.0494680603948896)-- (-0.5033495934959349,0.0494680603948896);
\draw (-0.5033495934959349,0.0494680603948896)-- (-0.5033495934959347,4.169468060394882);
\draw (-0.5033495934959347,4.169468060394882)-- (-4.623349593495928,4.169468060394884);
\draw (-4.623349593495928,4.169468060394884)-- (-4.623349593495929,0.0494680603948896);
\draw [dash pattern=on 3pt off 3pt] (3.595828106852492,4.158936120789778)-- (-0.5033495934959347,4.169468060394882);
\draw [dash pattern=on 3pt off 3pt] (3.595828106852492,4.158936120789778)-- (-0.5033495934959349,0.0494680603948896);
\draw (-4.623349593495928,4.169468060394884)-- (-0.5033495934959349,0.0494680603948896);
\draw (-0.5033495934959347,4.169468060394882)-- (-4.623349593495929,0.0494680603948896);
%\draw [dash pattern=on 3pt off 3pt] (3.595828106852492,4.158936120789778)-- (-4.623349593495929,0.0494680603948896);
\draw [dash pattern=on 3pt off 3pt] (-4.623349593495929,0.0494680603948896)-- (1.7070476190476158,-2.136998838559813);
\draw [dash pattern=on 3pt off 3pt] (-0.5033495934959347,4.169468060394882)-- (1.7070476190476158,-2.136998838559813);
%\draw [dash pattern=on 3pt off 3pt] (-0.5033495934959349,0.0494680603948896)-- (1.7070476190476158,-2.136998838559813);
\begin{scriptsize}
\draw [fill=qqqqff] (-4.623349593495929,0.0494680603948896) circle (2.5pt);
\draw[color=qqqqff] (-5,-0.5) node {$A$};
\draw [fill=qqqqff] (-0.5033495934959349,0.0494680603948896) circle (2.5pt);
\draw[color=qqqqff] (-0.2,-0.5) node {$B$};
\draw [fill=qqqqff] (-0.5033495934959347,4.169468060394882) circle (2.5pt);
\draw[color=qqqqff] (-0.2,4.6) node {$C$};
\draw [fill=qqqqff] (-4.623349593495928,4.169468060394884) circle (2.5pt);
\draw[color=qqqqff] (-5,4.6) node {$D$};
\draw [fill=qqqqff] (3.595828106852492,4.158936120789778) circle (2.5pt);
\draw[color=qqqqff] (4.1,4.6) node {$E$};
\draw [fill=qqqqff] (1.7070476190476158,-2.136998838559813) circle (2.5pt);
\draw[color=qqqqff] (2.1,-2.5) node {$E'$};
\end{scriptsize}
\end{tikzpicture}
\caption{Options for adding a fifth point to $E$ to the outside of a square.  Either option generates three distinct triangles.}
\label{fig: SquareAdding}
\end{subfigure}
%%%%%%%%%%%%%%%%%%%%%%%%%%%%%%%%%%%%%%%%%%%%%%%%%%%%%%%%%%%%%%
\caption{Options for adding a fifth point to a square.  Any choice except for the center of the square will result in a configuration with at least three distinct triangles.}
\end{figure}

\begin{proof}[Proof of case \ref{part: FourSides}: four congruent sides.]
There are two subcases: the four points either form a non-square rhombus or a square.

%\underline{Subcase \ref{part: FourSides}a: The Non-Square Rhombus}. 

If the four points form a non-square rhombus, then the argument presented in case \ref{part: TwoPairOpposite} for a non-rectangle parallelogram also applies to show that the addition of a fifth point anywhere generates a third triangle (see Figure \ref{fig: AddingTwoPairOpposite}).

If the four points form a square, we must show that the addition of a fifth point anywhere but the center results in a configuration determining at least three triangles.  If the fifth point is on the interior of the square but not in the center, then it creates a non-convex configuration (Figure \ref{fig: AddingSquareInside}).

%For the square, there are two possibilities. Either a point is added inside or outside. If we add a point inside, it must be in the center, where it determines exactly two triangles. Otherwise, the uneven distance between two of the sides necessarily creates at least three triangles.

If the fifth point $E$ is added outside the square, to avoid three distinct triangles, we must place it so that either $\tri EBC \cong \tri BCD$ or $\tri EBC \cong \tri EBA$ (see Figure \ref{fig: SquareAdding}).  If $\tri EBC \cong \tri BCD$, then $ECD$ are collinear, so there are at least three triangles.  If $\tri EBC \cong \tri EBA$, then we have a non-convex configuration, so there are at least three distinct triangles in this case also.

This shows that the addition of a fifth point to a square anywhere but the center generates at least three distinct triangles, and this completes the proof of case \ref{part: FourSides}.
\end{proof}

\begin{comment}
if it makes a 4-point subconfiguration that determines more than one new triangle, the 5-point total configuration will have at least three triangles. If we place the point $E$ such that $\tri EBC \cong \tri BCD$, then $E,C,D$ lie on a line. At least two new triangles appear: $\tri EBD$ and $\tri  EAD$. 
If we place the point $E$ such that $\tri FAB \cong \tri FBC$, $D,B,F$ are collinear. In addition to the new triangle $\tri FAB$, at least $\tri FAD$ is also new. Therefore, any point added outside of the square determines at least two new triangles. 
\end{comment}

\section{Proof of Theorem \ref{thm: DistinctTrianglesLowerBound}}\label{sec: TheoremProofDistinctTrianglesLowerBound}

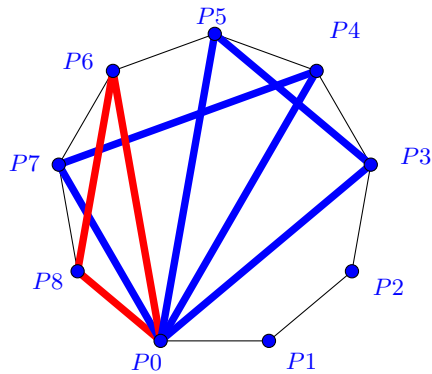
\begin{figure}[h]
\begin{tikzpicture}[line cap=round,line join=round,>=triangle 45,x=0.3cm,y=0.3cm]
\clip(-3.716,-10.63) rectangle (15.952,7.74);
\draw (3.22,-7.78)-- (8.021499765698216,-7.78);
\draw [line width=2.8pt,color=ffqqqq] (3.22,-7.78)-- (-0.4581622141501924,-4.693655442696365);
\draw (-1.2919338985318785,0.03489875264953923)-- (-0.4581622141501924,-4.693655442696365);
\draw (-1.2919338985318785,0.03489875264953923)-- (1.108815984317232,4.193119526009222);
\draw (5.620749882849115,5.835329164051487)-- (1.108815984317232,4.193119526009222);
\draw (5.620749882849115,5.835329164051487)-- (10.132683781380994,4.193119526009213);
\draw (12.533433664230095,0.03489875264952591)-- (10.132683781380994,4.193119526009213);
\draw (12.533433664230095,0.03489875264952591)-- (11.6996619798484,-4.693655442696375);
\draw (8.021499765698202,-7.78)-- (11.6996619798484,-4.693655442696375);
\draw [line width=2.8pt,color=qqqqff] (3.22,-7.78)-- (-1.2919338985318785,0.03489875264953923);
\draw [line width=2.8pt,color=qqqqff] (-1.2919338985318785,0.03489875264953923)-- (10.132683781380994,4.193119526009213);
\draw [line width=2.8pt,color=qqqqff] (10.132683781380994,4.193119526009213)-- (3.22,-7.78);
\draw [line width=2.8pt,color=ffqqqq] (-0.4581622141501924,-4.693655442696365)-- (1.108815984317232,4.193119526009222);
\draw [line width=2.8pt,color=ffqqqq] (1.108815984317232,4.193119526009222)-- (3.22,-7.78);
\draw [line width=2.8pt,color=qqqqff] (3.22,-7.78)-- (5.620749882849115,5.835329164051487);
\draw [line width=2.8pt,color=qqqqff] (5.620749882849115,5.835329164051487)-- (12.533433664230095,0.03489875264952591);
\draw [line width=2.8pt,color=qqqqff] (12.533433664230095,0.03489875264952591)-- (3.22,-7.78);
\begin{scriptsize}
\draw [fill=qqqqff] (3.22,-7.78) circle (2.5pt);
\draw[color=qqqqff] (2.598,-8.716) node {$P0$};
\draw [fill=qqqqff] (8.021499765698216,-7.78) circle (2.5pt);
\draw[color=qqqqff] (9.462,-8.65) node {$P1$};
\draw [fill=qqqqff] (-0.4581622141501924,-4.693655442696365) circle (2.5pt);
\draw[color=qqqqff] (-1.78,-5.13) node {$P8$};
\draw [fill=qqqqff] (-1.2919338985318785,0.03489875264953923) circle (2.5pt);
\draw[color=qqqqff] (-2.814,0.128) node {$P7$};
\draw [fill=qqqqff] (5.620749882849115,5.835329164051487) circle (2.5pt);
\draw[color=qqqqff] (5.48,6.64) node {$P5$};
\draw [fill=qqqqff] (1.108815984317232,4.193119526009222) circle (2.5pt);
\draw[color=qqqqff] (-0.438,4.616) node {$P6$};
\draw [fill=qqqqff] (10.132683781380994,4.193119526009213) circle (2.5pt);
\draw[color=qqqqff] (11.376,6.024) node {$P4$};
\draw [fill=qqqqff] (12.533433664230095,0.03489875264952591) circle (2.5pt);
\draw[color=qqqqff] (14.522,0.436) node {$P3$};
\draw [fill=qqqqff] (11.6996619798484,-4.693655442696375) circle (2.5pt);
\draw[color=qqqqff] (13.312,-5.372) node {$P2$};
\end{scriptsize}
\end{tikzpicture}
\caption{Illustrating the bijection described in the proof of Theorem \ref{thm: DistinctTrianglesLowerBound} with $n=9$.  Note that triangles $\tri P_0 P_4 P_7$ and $\tri P_0 P_3 P_5$ represent the same partition of $9$ $\left( \{4-0,\ 7-4,\ 9-7 \} = \{ 3-0,\ 5-3,\ 9-5\} = \{4,3,2\} \right)$.  Thus they are congruent; however, $\tri P_0 P_6 P_8$ represents a different partition $\left( \{6-0,\ 8-6,\ 9-8\} = \{6,2,1\} \right)$, so it is a different triangle.}
\label{fig: CountingRegularPolygon}
\end{figure}

\begin{proof}
We show that the vertices of a regular $n$-gon determine $[n^2/12]$ distinct triangles.  Conditional on Conjecture \ref{conj: RegularPolygons}, this completes the proof.  Label the vertices of a regular $n$-gon $\{P_0, \hdots, P_{n-1}\}$.  By the symmetry of the configuration, every congruence class of a triangle has a member with $P_0$ as a vertex, so when counting triangles we can just count triangles incident on $P_0$.  To form a triangle, we just have to pick two other vertices, $P_a$ and $P_b$, and we can assume $a<b$.  By symmetry, $\tri P_0 P_a P_b$ will be distinct from $\tri P_0 P_{a'} P_{b'}$ if and only if $\{a-0, b-a, n-b \}$ and $\{a'-0, b'-a', n-b'\}$ are not the same set (see Figure \ref{fig: CountingRegularPolygon}).  Thus there is a bijection between distinct triangles determined by the regular $n$-gon and ways to write $n$ as a sum of three positive integers.  Using a result from the theory of integer partitions (see \cite{Honsberger}), this quantity is equal to $[n^2/12]$, so this completes the proof.  A self-contained proof that this quantity is asymptotic to $n^2/12$ is also given in Appendix \ref{app: PartitionProof}.
\end{proof}

\section{Proof of Lemma \ref{lem: FourPointConfigs}} \label{sec: LemmaProof}

\begin{figure}[h]
\begin{tikzpicture}[line cap=round,line join=round,>=triangle 45,x=0.6cm,y=0.6cm]
\clip(0.7208284470274824,-0.43495488086704626) rectangle (6.3063607042351935,3.8866575276246715);
\draw (3.3822222222222225,3.671111111111111)-- (0.96,0.6);
\draw (0.96,0.6)-- (5.04,-0.32);
\draw (5.04,-0.32)-- (3.3822222222222225,3.671111111111111);
\draw (2.8422222222222224,1.171111111111111)-- (3.3822222222222225,3.671111111111111);
\draw (2.8422222222222224,1.171111111111111)-- (0.96,0.6);
\draw (2.8422222222222224,1.171111111111111)-- (5.04,-0.32);
\begin{scriptsize}
\draw [fill=qqqqff] (3.3822222222222225,3.671111111111111) circle (2.5pt);
\draw[color=qqqqff] (3.850271915245637,3.6769241029469124) node {$A$};
\draw [fill=qqqqff] (0.96,0.6) circle (2.5pt);
\draw[color=qqqqff] (0.8091372574181181,1.099410699670243) node {$B$};
\draw [fill=qqqqff] (5.04,-0.32) circle (2.5pt);
\draw[color=qqqqff] (5.517100711368886,-0.29697236463167853) node {$C$};
\draw [fill=qqqqff] (2.8422222222222224,1.171111111111111) circle (2.5pt);
\draw[color=qqqqff] (3.3204190529018223,1.237393215905611) node {$D$};
\end{scriptsize}
\end{tikzpicture}
\caption{Four points not in convex position; $\tri ABC$ and $\tri ABD$ are distinct.} \label{fig: NotConvex}
\end{figure}
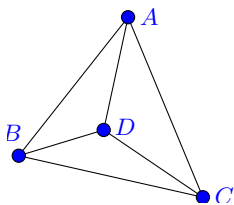

\begin{proof}[Proof of case \ref{part: NotConvex}: not in convex position.]
In this case, the four points form a triangle with one point in the interior (Figure \ref{fig: NotConvex}).  Triangle $\tri ABD$ is contained in $\tri ABC$, so they must be distinct.
\end{proof}

\begin{proof}[Proof of case \ref{part: Collinear}: three collinear points.]
Say point $C$ lies on $AB$ and $D$ does not (Figure \ref{fig: Collinear}).  Then $\tri ACD$ is contained in $\tri ABD$, so they are distinct.
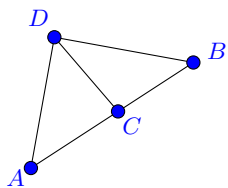
\begin{figure}[h]
\begin{tikzpicture}[line cap=round,line join=round,>=triangle 45,x=0.6cm,y=0.6cm]
\clip(2.983956112885607,1.5108379586786158) rectangle (8.432543095732711,5.544704110891731);
\draw (4.52,4.9)-- (4.,2.);
\draw (4.52,4.9)-- (7.6,4.34);
\draw (4.,2.)-- (7.6,4.34);
\draw (4.52,4.9)-- (5.930961383309082,3.2551248991509034);
\begin{scriptsize}
\draw [fill=qqqqff] (4.,2.) circle (2.5pt);
\draw[color=qqqqff] (3.6594534464315642,1.7784878455553155) node {$A$};
\draw [fill=qqqqff] (7.6,4.34) circle (2.5pt);
\draw[color=qqqqff] (8.126657510730768,4.601556890469076) node {$B$};
\draw [fill=qqqqff] (4.52,4.9) circle (2.5pt);
\draw[color=qqqqff] (4.169262754768136,5.334407771202896) node {$D$};
\draw [fill=qqqqff] (5.930961383309082,3.2551248991509034) circle (2.5pt);
\draw[color=qqqqff] (6.240363069885456,2.9383040220210135) node {$C$};
\end{scriptsize}
\end{tikzpicture}
\caption{Four points containting three collinear points; $\tri ACD$ and $\tri ABD$ are distinct.} \label{fig: Collinear}
\end{figure}
\end{proof}

\begin{proof}[Proof of case \ref{part: AllDistinct}: no congruent sides.]
Say the four points form quadrilateral $ABCD$ (Figure \ref{fig: AllDistinct}).  We have $\tri ABD \nc \tri CBD$ because $AB$, $AD$, $BC$, and $CD$ are all distinct.  We claim $\tri ABC$ is distinct from both of these.  Triangle $\tri ABC$ shares $AB$ with $\tri ABD$, and $BC \neq AD$, so if they are congruent then we must have $BC = BD$ and $AC = AD$.  This is impossible because then $\tri CBD$ and $\tri CAD$ would both be isoceles triangles with $CD$ as base, which is impossible unless one contains the other, which is not the case here.  Thus $\tri ABC \nc ABD$.  A similar argument shows that $\tri ABC \nc \tri CBD$, so we have three distinct triangles.
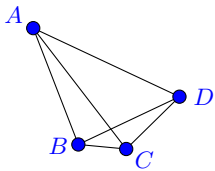
\begin{figure}[h]
\begin{tikzpicture}[line cap=round,line join=round,>=triangle 45,x=0.6cm,y=0.6cm]
\clip(1.7755024464664917,5.683830490709806) rectangle (7.09663710222943,9.539263385005116);
\draw (2.58,8.92)-- (3.58,6.34);
\draw (3.58,6.34)-- (4.64,6.24);
\draw (4.64,6.24)-- (5.82,7.4);
\draw (5.82,7.4)-- (2.58,8.92);
\draw (3.58,6.34)-- (5.82,7.4);
\draw (2.58,8.92)-- (4.64,6.24);
\begin{scriptsize}
\draw [fill=qqqqff] (2.58,8.92) circle (2.5pt);
\draw[color=qqqqff] (2.145114195010504,9.207887334586346) node {$A$};
\draw [fill=qqqqff] (3.58,6.34) circle (2.5pt);
\draw[color=qqqqff] (3.1328697299126063,6.3083468934221045) node {$B$};
\draw [fill=qqqqff] (4.64,6.24) circle (2.5pt);
\draw[color=qqqqff] (5.04465463617474,5.9897160757117485) node {$C$};
\draw [fill=qqqqff] (5.82,7.4) circle (2.5pt);
\draw[color=qqqqff] (6.363786221495612,7.4044369063457305) node {$D$};
\end{scriptsize}
\end{tikzpicture}
\caption{A quadrilateral with all distinct side lengths; $\tri ABC$, $\tri ABD$, and $\tri CBD$ are all distinct.} \label{fig: AllDistinct}
\end{figure}
\end{proof}

\begin{proof}[Proof of case \ref{part: OnePairAdjacent}: one pair of adjacent congruent sides.]
Let the points form quadrilateral $ABCD$ and suppose $AB = AD$ (Figure \ref{fig: OnePairAdjacent}).  Triangle $\tri ABD \nc \tri BCD$ because $\tri ABD$ is isoceles but $\tri BCD$ is not.  Also, by the same argument as in part \ref{part: AllDistinct}, we see that $\tri ABC$ is distinct from both of these, so there are at least three distinct triangles.
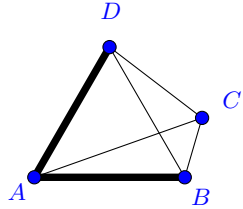
\begin{figure}[h]
\begin{tikzpicture}[line cap=round,line join=round,>=triangle 45,x=2.0cm,y=2.0cm]
\clip(2.760063002761122,6.637527966452467) rectangle (4.8115907614414395,8.123964246693903);
\draw [line width=2.8pt] (3.050025717307914,6.926490679780196)-- (4.050025717307914,6.926490679780196);
\draw [line width=2.8pt] (3.050025717307914,6.926490679780196)-- (3.550025717307914,7.792516083564634);
\draw (4.050025717307914,6.926490679780196)-- (4.165233579294158,7.321592893741037);
\draw (4.165233579294158,7.321592893741037)-- (3.550025717307914,7.792516083564634);
\draw (3.550025717307914,7.792516083564634)-- (4.050025717307914,6.926490679780196);
\draw (3.050025717307914,6.926490679780196)-- (4.165233579294158,7.321592893741037);
\begin{scriptsize}
\draw [fill=qqqqff] (3.050025717307914,6.926490679780196) circle (2.5pt);
\draw[color=qqqqff] (2.9345042852357364,6.826710765755922) node {$A$};
\draw [fill=qqqqff] (4.050025717307914,6.926490679780196) circle (2.5pt);
\draw[color=qqqqff] (4.1580501820295055,6.797227732098241) node {$B$};
\draw [fill=qqqqff] (3.550025717307914,7.792516083564634) circle (2.5pt);
\draw[color=qqqqff] (3.5585618309899885,8.033058226249386) node {$D$};
\draw [fill=qqqqff] (4.165233579294158,7.321592893741037) circle (2.5pt);
\draw[color=qqqqff] (4.364431417633274,7.443397553095758) node {$C$};
\end{scriptsize}
\end{tikzpicture}
\caption{Quadrilateral with one pair of adjacent congruent sides (shown in bold); $\tri ABD$, $\tri BCD$, and $\tri ABC$ are all distinct.}
\label{fig: OnePairAdjacent}
\end{figure}
\end{proof}

\begin{proof}[Proof of case \ref{part: OnePairOpposite}: one pair of opposite congruent sides.]
Suppose $AB = CD$ (Figure \ref{fig: OnePairOpposite}).  Triangle $\tri ABC \nc \tri DBC$ because they have two sides congruent to each other and the third is not.  We now claim that $\tri ACD$ is distinct from both of these.  Triangle $\tri ACD \nc \tri BCD$ by the same isoceles triangle argument from parts \ref{part: AllDistinct} and \ref{part: OnePairAdjacent}.  If $\tri ACD \cong \tri ABC$, then $BC$ must equal $AD$.  But that would force $AB$ to be parallel to $CD$, which would force $AC = BD$, a contradiction.  Thus there are at least three distinct triangles.
\begin{figure}[h]
\begin{tikzpicture}[line cap=round,line join=round,>=triangle 45,x=1.0cm,y=1.0cm]
\clip(2.3729822404398275,6.863389769213819) rectangle (5.918744745508566,9.432475177078006);
\draw [line width=2.8pt] (3.0352379835661467,7.288611626392931)-- (5.035237983566146,7.288611626392931);
\draw [line width=2.8pt] (3.,8.)-- (4.732050807568878,9.);
\draw (3.0352379835661467,7.288611626392931)-- (3.,8.);
\draw (5.035237983566146,7.288611626392931)-- (4.732050807568878,9.);
\draw (3.,8.)-- (5.035237983566146,7.288611626392931);
\draw (3.0352379835661467,7.288611626392931)-- (4.732050807568878,9.);
\begin{scriptsize}
\draw [fill=qqqqff] (3.0352379835661467,7.288611626392931) circle (2.5pt);
\draw[color=qqqqff] (2.8231029895862783,7.122421521081117) node {$A$};
\draw [fill=qqqqff] (5.035237983566146,7.288611626392931) circle (2.5pt);
\draw[color=qqqqff] (5.24781004395065,7.075710877301768) node {$B$};
\draw [fill=qqqqff] (3.,8.) circle (2.5pt);
\draw[color=qqqqff] (2.69146390257175,8.099098618285685) node {$C$};
\draw [fill=qqqqff] (4.732050807568878,9.) circle (2.5pt);
\draw[color=qqqqff] (5.0100104028921475,9.309328934386997) node {$D$};
\end{scriptsize}
\end{tikzpicture}
\caption{Quadrilateral with one pair of opposite congruent sides; $\tri ACD$, $\tri BCD$, and $\tri ABC$ are all distinct.}
\label{fig: OnePairOpposite}
\end{figure}
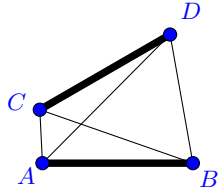
\end{proof}

\begin{proof}[Proof of case \ref{part: TwoPairAdjacent}: two pairs of adjacent congruent sides.]
Say $AB = AD$ and $BC = CD$ and assume without loss of generality that $AC > BD$ (Figure \ref{fig: TwoPairAdjacent}).  Triangle $\tri ABD \nc \tri BCD$ because $AB \neq BC$.  We claim that there is another triangle distinct from both of these.  First note that it is impossible to have both $AC = CD = BC$ and $BD = AD = AB$.  Because of this, the triangles $\tri ABD$, $\tri BCD$, and $\tri ACD$ are necessarily distinct, so there are at least three distinct triangles. 
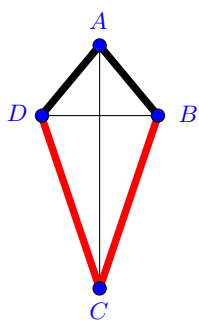
\begin{figure}[h]
\begin{tikzpicture}[line cap=round,line join=round,>=triangle 45,x=2.0cm,y=2.0cm]
\clip(2.409234753605918,7.415698386547161) rectangle (5.339616989199917,9.538909467426667);
\draw (3.812094711059942,9.266440070553228)-- (3.812094711059942,7.650251795787757);
\draw [line width=2.8pt] (3.812094711059942,9.266440070553228)-- (3.424189422119893,8.8);
\draw [line width=2.8pt] (3.812094711059942,9.266440070553228)-- (4.2,8.8);
\draw [line width=2.8pt,color=ffqqqq] (3.424189422119893,8.8)-- (3.812094711059942,7.650251795787757);
\draw [line width=2.8pt,color=ffqqqq] (4.2,8.8)-- (3.812094711059942,7.650251795787757);
\draw (3.424189422119893,8.8)-- (4.2,8.8);
\begin{scriptsize}
\draw [fill=qqqqff] (3.812094711059942,9.266440070553228) circle (2.5pt);
\draw[color=qqqqff] (3.802482355439233,9.430116833530361) node {$A$};
\draw [fill=qqqqff] (3.812094711059942,7.650251795787757) circle (2.5pt);
\draw[color=qqqqff] (3.8095012350454462,7.5039627010341465) node {$C$};
\draw [fill=qqqqff] (4.2,8.8) circle (2.5pt);
\draw[color=qqqqff] (4.402596561770459,8.812455428183597) node {$B$};
\draw [fill=qqqqff] (3.4290674320692798,8.799333632759735) circle (2.5pt);
\draw[color=qqqqff] (3.2655380655639252,8.815964867986702) node {$D$};
\end{scriptsize}
\end{tikzpicture}
\caption{Quadrilateral with two pairs of adjacent congruent sides.  Independently of the lengths of $AC$ and $BD$, the triangles $\tri ABD$, $\tri BCD$, and $\tri ACD$ are all distinct.} \label{fig: TwoPairAdjacent}
\end{figure}
\end{proof}

\begin{proof}[Proof of case \ref{part: ThreeSides}: three congruent sides.]
Say $AD = AB = BC$ (Figure \ref{fig: ThreeSides}).  Triangle $\tri ABC \nc \tri ADC$ because they have two sides congruent with each other and one side not congruent, thus there are at least two distinct triangles.
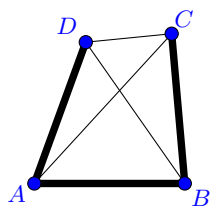
\begin{figure}[h]
\begin{tikzpicture}[line cap=round,line join=round,>=triangle 45,x=1.0cm,y=1.0cm]
\clip(2.3598042939791646,7.565937504518155) rectangle (6.260143049554784,10.391931453168779);
\draw [line width=2.8pt] (3.0514622375744236,8.02985035209082)-- (5.051462237574423,8.02985035209082);
\draw [line width=2.8pt] (3.0514622375744236,8.02985035209082)-- (3.735502524225761,9.909235593662636);
\draw [line width=2.8pt] (4.877150752079107,10.02223974827431)-- (5.051462237574423,8.02985035209082);
\draw (3.735502524225761,9.909235593662636)-- (4.877150752079107,10.02223974827431);
\draw (3.735502524225761,9.909235593662636)-- (5.051462237574423,8.02985035209082);
\draw (4.877150752079107,10.02223974827431)-- (3.0514622375744236,8.02985035209082);
\begin{scriptsize}
\draw [fill=qqqqff] (3.0514622375744236,8.02985035209082) circle (2.5pt);
\draw[color=qqqqff] (2.8175686030167824,7.8929120109736) node {$A$};
\draw [fill=qqqqff] (5.051462237574423,8.02985035209082) circle (2.5pt);
\draw[color=qqqqff] (5.265206337054655,7.841530302816316) node {$B$};
\draw [fill=qqqqff] (3.735502524225761,9.909235593662636) circle (2.5pt);
\draw[color=qqqqff] (3.4855308090614687,10.121009719248553) node {$D$};
\draw [fill=qqqqff] (4.877150752079107,10.02223974827431) circle (2.5pt);
\draw[color=qqqqff] (5.0409952469137815,10.219102071185187) node {$C$};
\end{scriptsize}
\end{tikzpicture}
\caption{Quadrilateral with three congruent sides; $\tri ABC$ and $\tri ADC$ are distinct.} \label{fig: ThreeSides}
\end{figure}
\end{proof}

\appendix

\section{Number of distinct triangles determined by a regular $n$-gon}
\label{app: PartitionProof}

We give a self-contained proof that the number of distinct triangles determined by a regular $n$-gon is asymptotic to $n^2/12$.  In the proof of Theorem \ref{thm: DistinctTrianglesLowerBound}, we establish that this is equal to the number of ways to write $n$ as a sum of three positive integers.  Denote this quantity by $p(n,3)$.  Since the order of a partition doesn't matter, we view this quantity as the number of ways to pick two elements $k < l$ from $\{1, \hdots, n\}$ such that $k \geq l-k \geq n-l > 0$.  Note that $k$ can be any of the elements $\ceiling{n/3}, \hdots, n-2$.  Once $k$ is chosen, $l$ can be any of the elements $k+\ceiling{(n-k)/2}, \hdots, \min(2k,\ n-1)$.  Note $2k$ is the minimum when $k \leq \floor{n/2}$, and $n-1$ is the minimum otherwise.  Thus the number of choices is given by
\begin{align}
p(n,3) \ &= \ \sum_{k = \ceiling{n/3}}^{\floor{n/2}} \ \sum_{l = k + \ceiling{(n-k)/2}}^{2k} 1 + \sum_{k = \floor{n/2}+1}^{n-2} \ \sum_{l = k + \ceiling{(n-k)/2}}^{n-1} 1 \nonumber \\
&= \ \sum_{k=n/3}^{n/2} \ \sum_{l=k+(n-k)/2}^{2k} 1 + \sum_{k=(n+2)/2}^{n-2} \ \sum_{l=k+(n-k)/2}^{n-1} 1 \ + \ O(n) \nonumber \\
\ &= \ \sum_{k=n/3}^{n/2} (3k/2 - n/2 + 1) + \sum_{k=(n+2)/2}^{n-2} (n/2 - k/2) \ + \ O(n) \nonumber \\
\ &= \ \frac{3}{4} \left( \frac{n^2}{4} - \frac{n^2}{9} \right) - \frac{n^2}{12} + \frac{n^2}{4} - \frac{1}{4} \left( n^2 - \frac{n^2}{4} \right) \ + \ O(n)   \nonumber \\
\ &= \ \frac{n^2}{12} \ + \ O(n), 
\end{align}
and this completes the proof. 
\hfill $\Box$
%%%%%%%%%%%%%%%%%%%%%%%%%%%%%%%%%%%%%%%%%%%%%%%%%%%%%%%%%%%%%%%%%
%%%%%%%%%%%%%%%%%%%%%%%%%%%%%%%%%%%%%%%%%%%%%%%%%%%%%%%%%%%%%%%%%

% note the form of the bibliography. start with a bibitem.
% you refer to an entry in the bibliography by using the label in braces
% the computer prints the label in brackets. often I have both the same, but
% you don't have to.

\bigskip

\end{document}